\documentclass[12pt]{article}
\usepackage{amsmath}
\usepackage{amssymb}
\usepackage{array}
\usepackage{geometry}
\usepackage{graphicx}
\usepackage{enumerate}
\usepackage[small, compact] {titlesec}
\numberwithin{equation}{section}
\begin{document}
\vskip7cm\noindent
\begin{center}
{\bf Evaluation of Matrix-variate Gamma and Beta Integrals}\\
\vskip.2cm{\bf as Multiple Integrals and Kober Fractional Integral Operators}\\
\vskip.2cm{\bf in the Complex Matrix Variate Case}\\
\vskip.3cm{A.M. Mathai}\\
 \vskip.2cm{Director, Centre for
Mathematical Sciences India}\\
 \vskip.1cm{[Arunapuram
P.O., Palai, Kerala-686574, Kerala, India]}\\
directorcms458@gmail.com , Phone:91+9495427558\\
\vskip.2cm{and}\\
 \vskip.2cm{Emeritus Professor of
Mathematics and Statistics, McGill University Canada;\\
mathai@math.mcgill.ca}\\
 \vskip.1cm{[805 Sherbrooke
Street West, Montreal, Quebec, Canada, H3A2K6]}
\end{center}

\vskip.5cm\noindent{\bf Abstract} \vskip.3cm Explicit evaluations of 
matrix-variate gamma and beta integrals in the complex domain by using
conventional procedures is extremely difficult. Such an evaluation will
reveal the structure of these matrix-variate integrals. In this article, explicit evaluations
of matrix-variate gamma and beta integrals in the complex domain for the order
of the matrix $p=1,2$ are given. Then fractional integral operators of the
Kober type are given for some specific cases of the arbitrary function.
A formal definition of
fractional integrals in the complex matrix-variate case was given
by the author earlier as the M-convolution of products and ratios,
where Kober operators become a special class
of fractional integral operators.



\vskip.5cm\noindent{\bf 1.\hskip.3cm Introduction} \vskip.3cm

There is intensified activity in the area of fractional calculus in recent years due to its many applications in engineering, physical, biological and social sciences. Solutions coming out of fractional differential equations are seen to describe real-life and experimental situations much better compared to the solutions of conventional differential equations. The reason seems to be that fractional derivatives are certain types of integrals and thereby describe global activity whereas conventional derivatives describe local activity. Fractional calculus in the scalar and matrix-variate case in the complex domain was not  available in the literature but recently given in [8],[9]. Fractional calculus for the real scalar variable case is available in many books and articles, see for example [1],[2],[4],[11],[14]. Solutions of fractional differential equations in the real scalar case is available in many books and journals, see for example, [3]. Fractional integral operators in the scalar and real matrix-variate cases may be seen from [7],[12]. A pathway extension may be seen from [6],[10].
\vskip.2cm
This author has given a formal definition of fractional integrals in the real and complex matrix-variate case recently as M-convolutions of products and ratios, see [8], [9]. Let $\tilde{X_1}$ and $\tilde{X_2}$ be $p\times p$ hermitian positive definite matrices and let $\tilde{U_1}=\tilde{X_2}^{\frac{1}{2}}\tilde{X_1}\tilde{X_2}^{\frac{1}{2}}$ and $\tilde{U_2}=\tilde{X_2}^{\frac{1}{2}}\tilde{X_1}^{-1}\tilde{X_2}^{\frac{1}{2}}$. Then $\tilde{U_1}$ is called the product and $\tilde{U_2}$ is called the ratio. M-convolutions correspond to Mellin convolutions of products and ratios in the scalar variable cases. Product will lead to type-2 or right-sided fractional integrals and the ratio will lead to type-1 or left-sided fractional integrals. Fractional integrals of the Kober type or Kober operators are denoted by $K_{2,\tilde{U}}^{-\alpha}f$ and $K_{1,\tilde{U}}^{-\alpha}f$ respectively, where $\alpha$ indicates the order of the integrals. As per the new definition introduced, these Kober operators have the following representations:
$$K_{2,\tilde{U}}^{-\alpha}f=\frac{|{\rm det}(\tilde{U})|^{\beta}}{\tilde{\Gamma_p}(\alpha)}\int_{\tilde{V}>\tilde{U}}|{\rm det}(\tilde{V})|^{-\beta-\alpha}|{\rm det}(\tilde{V}-\tilde{U})|^{\alpha-p}f(\tilde{V}){\rm d}\tilde{V}.\eqno(1.1)$$
$$K_{1,\tilde{U}}^{-\alpha}f=\frac{|{\rm det}(\tilde{U})|^{-\alpha-\beta}}{\tilde{\Gamma_p}(\alpha)}\int_{\tilde{V}<\tilde{U}}|{\rm det}(\tilde{V})|^{\beta}|{\rm det}(\tilde{U}-\tilde{V})|^{\alpha-p}f(\tilde{V}){\rm d}\tilde{V}.\eqno(1.2)$$ In this paper we will examine (1.1) and (1.2) for many cases of the arbitrary function $f$ after evaluating matrix-variate gamma and beta integrals explicitly for the cases $p=1,2$ in the complex domain. Explicit evaluations of matrix-variate integrals are difficult and the evaluations given in this article reveal the structure of these matrix-variate integrals in the real and complex cases.

\vskip.2cm We will use the following standard notations.  All matrices appearing are $p\times p$ with elements in the
complex domain unless otherwise stated. Determinant of $(\cdot)$ will be denoted by ${\rm det}(\cdot)$ and $|{\rm det}(\cdot)|$ will denote
the absolute value of the determinant of $(\cdot)$. Matrices in the complex domain will be written with a tilde,  $\tilde{X}$.
Constant matrices will not be written with a tilde whether in the
real or complex domain. ${\rm tr}(X)$ is the trace of $X$, $({\rm
d}\tilde{X})=({\rm d}\tilde{x}_{ij})$ is the matrix of differentials
${\rm d}\tilde{x}_{ij}$'s. Let $\tilde{X}=X_1+iX_2$ where $X_1$ and
$X_2$ are real $m\times n$ matrices and $i=\sqrt{-1}$. Then ${\rm
d}\tilde{X}={\rm d}X_1\wedge{\rm d}X_2$ where
$${\rm d}X_1=\prod_{i=1}^m\prod_{j=1}^n\wedge{\rm d}x_{ij1}\mbox{
and }{\rm d}X_2=\prod_{i=1}^m\prod_{j=1}^n\wedge{\rm
d}x_{ij2}\nonumber
$$where $x_{ij1}$ and $x_{ij2}$ are the $(i,j)-th$ elements in $X_1$
and $X_2$ respectively, and $\wedge$ denotes the wedge product. For
any $p\times p$ matrix $B=B_1+iB_2$ in the complex domain, the
determinant will be a complex number of the form ${\rm det}(B)=b+ic$
where $b$ and $c$ are real scalar quantities. Then the absolute
value of the determinant will be of the form $|{\rm
det}(B)|=[(b+ic)(b-ic)]^{\frac{1}{2}}=[b^2+c^2]^{\frac{1}{2}}$. Note
that the conjugate of $B_1+iB_2$ is $B_1-iB_2$.

\vskip.2cm We need a few basic results on Jacobians of matrix transformations
in the complex domain. These results, other results  and  properties may be seen from
[5]. The results that we need will be listed here as lemmas.

\vskip.3cm\noindent{\bf Lemma 1.1.}\hskip.3cm{\it Let $\tilde{X}$
and $\tilde{Y}$ be $m\times n$ matrices in the complex domain. Let
$A$ be $m\times m$ and $B$ be $n\times n$ nonsingular constant
matrices in the sense of free of the elements in $\tilde{X}$ and
$\tilde{Y}$. Let $C$ be a constant $m\times n$ matrix. Then
$$\tilde{Y}=A\tilde{X}B+C,{\rm det}(A)\ne 0,{\rm det}(B)\ne
0\Rightarrow {\rm d}\tilde{Y}=|{\rm det}(AA^{*})|^n|{\rm
det}(BB^{*})|^m{\rm d}\tilde{X},\eqno(1.3)
$$where $A^{*}$ and $B^{*}$ denote the conjugate transposes of $A$
and $B$ respectively.}

\vskip.3cm When $A=A^{*}$, where $A^{*}$ denotes the conjugate transpose of $A$, then it is called a
hermitian matrix. The next result is about the transformation of a
hermitian matrix to a hermitian matrix.

\vskip.3cm\noindent{\bf Lemma 1.2.}\hskip.3cm{\it Let $\tilde{X}$
and $\tilde{Y}$ be $p\times p$ hermitian matrices and let $A$ be a
nonsingular constant matrix. Then
$$\tilde{Y}=A\tilde{X}A^{*}\Rightarrow {\rm d}\tilde{Y}=\begin{cases}|{\rm
det}(A)|^{2p}{\rm d}\tilde{X}\\ |{\rm det}(AA^{*})|^p{\rm
d}\tilde{X}\end{cases}\eqno(1.4)
$$}
\vskip.2cm The next result is on a decomposition of the hermitian
positive definite matrix $\tilde{X}=\tilde{X}^{*}>O$.

\vskip.3cm\noindent{\bf Lemma 1.3.}\hskip.3cm{\it Let $\tilde{X}$ be
a $p\times p$ hermitian positive definite matrix. Let $\tilde{T}$ be
a $p\times p$ lower triangular matrix with diagonal elements
$t_{jj}$'s being real and positive. Consider the unique
representation $\tilde{X}=\tilde{T}\tilde{T}^{*}$. Then
$$\tilde{X}=\tilde{T}\tilde{T}^{*}\Rightarrow {\rm
d}\tilde{X}=2^p\{\prod_{j=1}^pt_{jj}^{2(p-j)+1}\}{\rm
d}\tilde{T}.\eqno(1.5)
$$}
\vskip.3cm Next we define a complex matrix variate gamma function,
denoted by $\tilde{\Gamma}_p(\alpha)$ and defined as
$$\tilde{\Gamma}_p(\alpha)=\pi^{\frac{p(p-1)}{2}}\Gamma(\alpha)\Gamma(\alpha-1)...\Gamma(\alpha-p+1),~\Re(\alpha)>p-1,\eqno(1.6)
$$with the following integral
representation:
$$\tilde{\Gamma}_p(\alpha)=\int_{\tilde{Y}>O}|{\rm
det}(\tilde{Y})|^{\alpha-p}{\rm e}^{-{\rm tr}(\tilde{Y})}{\rm
d}\tilde{Y}.\eqno(1.7)
$$By using Lemma 1.3 and (1.6) we can evaluate the integral in (1.7). Then a matrix-variate gamma density, in the complex domain, can be defined  as follows:
$$f(\tilde{X})={{|{\rm
det}(B)|^{\alpha}}\over{\tilde{\Gamma}_p(\alpha)}}|{\rm
det}(\tilde{X})|^{\alpha-p}{\rm e}^{-{\rm
tr}(B\tilde{X})},\tilde{X}=\tilde{X}^{*}>O,~\Re(\alpha)>p-1\eqno(1.8)
$$and $f(\tilde{X})=0$ elsewhere, where $B=B^{*}>O$ is a constant
hermitian positive definite matrix.

\vskip.3cm\noindent{\bf Lemma 1.4}.\hskip.3cm{\it Let $\tilde{X}$ be
a nonsingular matrix and let $\tilde{Y}=\tilde{X}^{-1}$. Then
$$\tilde{Y}=\tilde{X}^{-1}\Rightarrow {\rm d}\tilde{Y}=\begin{cases}|{\rm
det}(\tilde{X}\tilde{X}^{*})|^{-2p}{\rm d}\tilde{X}\mbox{ for a
general $\tilde{X}$}\\ |{\rm
det}(\tilde{X}\tilde{X}^{*})|^{-p}\mbox{ for
$\tilde{X}=\tilde{X}^{*}$ or
$\tilde{X}=-\tilde{X}^{*}$}.\end{cases}\eqno(1.9)
$$}We need complex matrix variate beta function and its integral
representations. The complex matrix variate beta function will be
denoted and defined as follows:
\begin{align}
\tilde{B}_p(\alpha,\beta)&=\frac{\tilde{\Gamma}_p(\alpha)\tilde{\Gamma}_p(\beta)}{\tilde{\Gamma}_p(\alpha+\beta)},
\Re(\alpha)>p-1,\Re(\beta)>p-1&(1.10)\nonumber\\&=\int_{O<\tilde{X}<I}|{\rm
det}(\tilde{X})|^{\alpha-p}|{\rm det}(I-\tilde{X})|^{\beta-p}{\rm
d}\tilde{X}&(1.11)\nonumber\\ &=\int_{\tilde{U}>0}|{\rm
det}(\tilde{U})|^{\alpha-p}|{\rm
det}(I+\tilde{U})|^{-(\alpha+\beta)}{\rm d}\tilde{U}&(1.12)\nonumber
\end{align} for $\Re(\alpha)>p-1,\Re(\beta)>p-1$
where, in general, $\int_{A<\tilde{X}<B}f(\tilde{X}){\rm
d}\tilde{X}$ will mean the integral of a real-valued scalar function
$f(\tilde{X})$ of complex matrix argument $\tilde{X}$ and the
integral is taken over all $\tilde{X}$ such that
$A=A^{*}>O,B=B^{*}>O, \tilde{X}=\tilde{X}^{*}>O,
\tilde{X}-A>O,B-\tilde{X}>O$, where $A$ and $B$ are constant
matrices.

\vskip.3cm\noindent{\bf 2.\hskip.3cm Explicit Evaluations of Gamma and Beta Integrals in the
Matrix Variate Case in the Complex Domain}

\vskip.3cm We will try to evaluate explicitly matrix-variate gamma and beta integrals in the complex domain by using conventional procedures. This will reveal the structure of the integrals. When the real $p\times p$ matrix $X$ is symmetric and positive definite then for $p=1$ it reduces to a real scalar positive variable $x>0$. If the matrix $\tilde{X}$ is in the complex domain and if it is hermitian positive definite then also for $p=1$ it reduces to a real scalar positive variable $x>0$. Hence for $p=1$ the matrix-variate gamma and beta integrals in the real positive definite or hermitian positive definite cases reduce to the ordinary gamma and beta integrals in the real scalar case.
\vskip.3cm\noindent{\bf 2.1.\hskip.3cm Matrix-Variate Gamma in the Real Case, $p=2$}
\vskip.3cm
For $p=2$
$$\Gamma_2(\alpha)=\pi^{\frac{1}{2}}\Gamma(\alpha)\Gamma(\alpha-\frac{1}{2})=\int_{X>O}|X|^{\alpha-\frac{3}{2}}{\rm e}^{-{\rm tr}(X)}{\rm d}X.
$$Let
$$X=\left[\begin{matrix}x_1&x_2\\
x_2&x_3\end{matrix}\right]>O\Rightarrow x_1>0,x_3>0,x_1x_3-x_2^2>0.$$Hence we have to evaluate a triple integral over $x_1,x_2,x_3$ subject to the conditions $x_1>0,x_3>0,x_1x_3-x_2^2>0$. The integral, denoted by $\Gamma_2$ is the following:

$$\Gamma_2=\int\int\int[x_1x_3-x_2^2]^{\alpha-\frac{3}{2}}{\rm e}^{-(x_1+x_3)}{\rm d}x_1\wedge{\rm d}x_2\wedge{\rm d}x_2.$$Let $x_2=\sqrt{x_1x_3}~r$ for fixed $x_1,x_3$, then the Jacobian is $\sqrt{x_1x_3}$. The integral becomes
$$\Gamma_2=\int\int\int(x_1x_3)^{\alpha-\frac{3}{2}}(1-r^2)^{\alpha-\frac{3}{2}}{\rm e}^{-(x_1+x_3)}\sqrt{x_1x_3}{\rm d}x_1\wedge{\rm d}x_3\wedge{\rm d}r
$$for $-1<r<1,x_1>0,x_3>0.$ Integral over $r$ gives
\begin{align}
\int_{-1}^{1}(1-r^2)^{\alpha-\frac{3}{2}}{\rm d}r&=2\int_0^1(1-r^2)^{\alpha-\frac{3}{2}}{\rm d}r,~u=r^2\nonumber\\
&=\frac{\Gamma(\frac{1}{2})\Gamma(\alpha-\frac{1}{2})}{\Gamma(\alpha)}=\sqrt{\pi}\frac{\Gamma(\alpha-\frac{1}{2})}
{\Gamma(\alpha)}.\nonumber
\end{align} But
$$\int_{x_1>0}x_1^{\alpha-1}{\rm e}^{-x_1}{\rm d}x_1=\int_{x_3>0}x_3^{\alpha-1}{\rm e}^{-x_3}{\rm d}x_3=\Gamma(\alpha).
$$Hence the integral is

\begin{align}
\Gamma_2&=[\Gamma(\alpha)]^2\pi^{\frac{1}{2}}\frac{\Gamma(\alpha-\frac{1}{2})}{\Gamma(\alpha)}\nonumber\\
&=\pi^{\frac{1}{2}}\Gamma(\alpha)\Gamma(\alpha-\frac{1}{2}),~\Re(\alpha)>\frac{1}{2}\nonumber\\
&=\Gamma_2(\alpha)\nonumber
\end{align} and hence the result is verified.

\vskip.3cm\noindent{\bf 2.2.\hskip.3cm Matrix-Variate Gamma Integral in the Complex Case, $p=2$}
\vskip.3cm
From our notation
$$\tilde{\Gamma_2}(\alpha)=\pi \Gamma(\alpha)\Gamma(\alpha-1),~\Re(\alpha)>1.
$$Consider the matrix
$$\tilde{X}=\left[\begin{matrix}
x_1&x_2+iy_2\\
x_2-iy_2&x_3
\end{matrix}\right]=\tilde{X}^{*}
$$where * denotes the conjugate transpose. When $\tilde{X}>O$ or hermitian positive definite then we have $x_1>0,x_3>0,x_1x_3-(x_2^2+y_2^2)>0$. Hence integration is to be done under these conditions. Let the integral be denoted by $\tilde{\Gamma_2}$. Then
\begin{align}
\tilde{\Gamma_2}&=\int_{\tilde{X}>O}|{\rm det}(\tilde{X})|^{\alpha-2}{\rm e}^{-{\rm tr}(\tilde{X})}{\rm d}\tilde{X}\nonumber\\
&=\int\int\int\int [x_1x_3-(x_2^2+y_2^2)]^{\alpha-2}{\rm e}^{-(x_1+x_3)}{\rm d}x_1\wedge{\rm d}x_3\wedge{\rm d}x_2\wedge{\rm d}y_2.\nonumber
\end{align} Let $x_2=\sqrt{x_1x_3}~r\cos\theta,y_2=\sqrt{x_1x_3}~r\sin\theta$ then the Jacobian is $x_1x_3~r$ and the quadruple integral becomes
$$\tilde{\Gamma_2}=\int_{\theta=0}^{2\pi}\int_{r=0}^1\int_{x_1>0}\int_{x_3>0}(x_1x_3)^{\alpha-1}r(1-r^2)^{\alpha-2}{\rm e}^{-(x_1+x_3)}{\rm d}x_1\wedge{\rm d}x_3\wedge{\rm d}r\wedge{\rm d}\theta.
$$Integral over $\theta$ gives $2\pi$. Put $u=r^2$. Integral over $r$ gives
$$2\int_0^1r(1-r^2)^{\alpha-2}{\rm d}r=\int_0^1(1-u)^{\alpha-2}{\rm d}u=\frac{\Gamma(1)\Gamma(\alpha-1)}{\Gamma(\alpha)}.
$$Integrals over $x_1$ and $x_3$ give one $\Gamma(\alpha)$ each. Then the quadruple integral gives
$$\pi[\Gamma(\alpha)]^2\frac{\Gamma(1)\Gamma(\alpha-1)}{\Gamma(\alpha)}=\pi~\Gamma(\alpha)\Gamma(\alpha-1)
=\tilde{\Gamma_2}(\alpha)
$$for $\Re(\alpha)>1$ which verifies the result.

\vskip.3cm\noindent{\bf 2.3.\hskip.3cm Matrix-Variate Beta Integral in the Real Case, $p=2$}
\vskip.3cm The general definition of real matrix-variate beta function and an integral representation are the following:
$$B_p(\alpha,\beta)=\frac{\Gamma_p(\alpha)\Gamma_p(\beta)}{\Gamma_p(\alpha+\beta)}=\int_{O}^{I}|X|^{\alpha-\frac{p+1}{2}}|I-X|^{\beta-\frac{p+1}{2}}{\rm d}X
$$for $X>O,\Re(\alpha)>\frac{p-1}{2},\Re(\beta)>\frac{p-1}{2}$. As mentioned earlier, for $p=1$ the real and complex matrix-variate cases coincide with the real scalar variable case. Hence we look into the case $p=2$.
\begin{align}
B_2(\alpha,\beta)&=\frac{\Gamma_2(\alpha)\Gamma_2(\beta)}{\Gamma_2(\alpha+\beta)}\nonumber\\
&=\pi^{\frac{1}{2}}\frac{\Gamma(\alpha)\Gamma(\alpha-\frac{1}{2})\Gamma(\beta)\Gamma(\beta-\frac{1}{2})}
{\Gamma(\alpha+\beta)\Gamma(\alpha+\beta-\frac{1}{2})}\nonumber
\end{align}for $\Re(\alpha)>\frac{1}{2},\Re(\beta)>\frac{1}{2}.$ Let
$$X=\left[\begin{matrix}x_1&x_2\\
x_2&x_3
\end{matrix}\right]\mbox{ then }I-X=\left[\begin{matrix}1-x_1&-x_2\\
-x_2&1-x_3
\end{matrix}\right]$$and the integral representation becomes
$$\int_{O}^{I}|X|^{\alpha-\frac{3}{2}}|I-X|^{\beta-\frac{3}{2}}{\rm d}X=\int_{O}^{I}[x_1x_3-x_2^2]^{\alpha-\frac{3}{2}}[(1-x_1)(1-x_3)-x_2^2]^{\beta-\frac{3}{2}}{\rm d}x_1\wedge{\rm d}x_3\wedge{\rm d}x_2.
$$But
$$(x_1x_3-x_2^2)^{\alpha-\frac{3}{2}}=x_3^{\alpha-\frac{3}{2}}[x_1-\frac{x_2^2}{x_3}]^{\alpha-\frac{3}{2}}
$$which means that $x_1>\frac{x_2^2}{x_3}$. Now,
$$[(1-x_1)(1-x_3)-x_2^2]^{\beta-\frac{3}{2}}=(1-x_3)^{\beta-\frac{3}{2}}[1-x_1-\frac{x_2^2}{1-x_3}]^{\beta-\frac{3}{3}}
$$which means that $x_1<1-\frac{x_2^2}{1-x_3}$. That is,
$$\frac{x_2^2}{x_3}<x_1<1-\frac{x_2^2}{1-x_3}
$$which means that $0<u<b$ where $u=x_1-\frac{x_2^2}{x_3},~b=1-\frac{x_2^2}{x_3(1-x_3)}$. Then the integrand will reduce to the following factors:
\begin{align}
x_3^{\alpha-\frac{3}{2}}&(1-x_3)^{\beta-\frac{3}{2}}u^{\alpha-\frac{3}{2}}b^{\beta-\frac{3}{2}}
[1-\frac{u}{b}]^{\beta-\frac{3}{2}}\nonumber\\
&=x_3^{\alpha-\frac{3}{2}}(1-x_3)^{\beta-\frac{3}{2}}v^{\alpha-\frac{3}{2}}(1-v)^{\beta-\frac{3}{2}}b^{\alpha+\beta-2}\nonumber
\end{align}where $v=\frac{u}{b}$. Put $z=\frac{x_2}{\sqrt{x_3(1-x_3)}}$ for fixed $x_3$. The integral over $z$ gives
\begin{align}
\int_z(1-z^2)^{\alpha+\beta-2}{\rm d}z&=2\int_{z>0}(1-z^2)^{\alpha+\beta-2}{\rm d}z\nonumber\\
&=\int_{w>0}w^{\frac{1}{2}-1}(1-w)^{\alpha+\beta-2}{\rm d}w,w=z^2\nonumber\\
&=\frac{\Gamma(\frac{1}{2})\Gamma(\alpha+\beta-1)}{\Gamma(\alpha+\beta-\frac{1}{2})},\Re(\alpha+\beta)>1.\nonumber
\end{align}Now the integrals over $x_3, v$ and $w$ give the following:
\begin{align}
\int_0^1x_3^{\alpha-1}(1-x_3)^{\beta-1}{\rm d}x_3&=\frac{\Gamma(\alpha)\Gamma(\beta)}{\Gamma(\alpha+\beta)},\Re(\alpha)>0,\Re(\beta)>0\nonumber\\
\int_0^1v^{\alpha-\frac{3}{2}}(1-v)^{\beta-\frac{3}{2}}{\rm d}v&=\frac{\Gamma(\alpha-\frac{1}{2})\Gamma(\beta-\frac{1}{2})}{\Gamma(\alpha+\beta-1)},\Re(\alpha)>\frac{1}{2},
\Re(\beta)>\frac{1}{2}\nonumber\\
2\int_{z>0}(1-z^2)^{\alpha+\beta-2}{\rm d}z&=\frac{\Gamma(\frac{1}{2})\Gamma(\alpha+\beta-1)}{\Gamma(\alpha+\beta-\frac{1}{2})}.\nonumber
\end{align}Hence the total integral, by taking the product, is
$$\pi^{\frac{1}{2}}\frac{\Gamma(\alpha)\Gamma(\alpha-\frac{1}{2})\Gamma(\beta)\Gamma(\beta-\frac{1}{2})}{\Gamma(\alpha+\beta)
\Gamma(\alpha+\beta-\frac{1}{2})}=\frac{\Gamma_2(\alpha)\Gamma_2(\beta)}{\Gamma_2(\alpha+\beta)}=B_2(\alpha,\beta).
$$Hence the result is verified for the real case. Explicit evaluations for $p\ge 3$ will be difficult and it is simpler to use matrix methods directly starting from $p\ge 2$.

\vskip.3cm\noindent{\bf 2.4.\hskip.3cm Matrix-variate Beta Integral in the Complex Case: $p=2$}
\vskip.3cm Let
$$\tilde{X}=\left[\begin{matrix}
x_1&x_2+iy_2\\
x_2-iy_2&x_3
\end{matrix}\right],~x_1>0,x_3>0$$and
\begin{align}
|{\rm det}(\tilde{X})|^{\alpha-2}&=[x_1x_3-(x_2^2+y_2^2)]^{\alpha-2}\nonumber\\
|{\rm det}(I-\tilde{X})|^{\beta-2}&=[(1-x_1)(1-x_3)-(x_2^2+y_2^2)]^{\beta-2}.\nonumber
\end{align}Steps parallel to the ones in the real case will go through and we have the factors
$$x_3^{\alpha-2}(1-x_3)^{\beta-2}v^{\alpha-2}(1-v)^{\beta-2}b^{\alpha+\beta-3},b=1-\frac{(x_2^2+y_2^2)}{x_3(1-x_3)}.
$$Put $z_1=\frac{x_2}{\sqrt{x_3(1-x_3)}},~z_2=\frac{y_2}{\sqrt{x_3(1-x_3)}}$ for fixed $x_3$. Then $b=1-(z_1^2+z_2^2)$. Put $z_1=r~\cos\theta,~z_2=r~\sin\theta,~0\le r\le 1, ~0\le\theta\le 2\pi$. The Jacobian is $r$.
\begin{align}
\int_{z_1}\int_{z_2}b^{\alpha+\beta-3}{\rm d}z_1\wedge{\rm d}z_2&=\int_{0}^{2\pi}\int_{r=0}^1r(1-r^2)^{\alpha+\beta-3}{\rm d}r\wedge{\rm d}\theta\nonumber\\
&=\pi\int_0^12r(1-r^2)^{\alpha+\beta-3}{\rm d}r\nonumber\\
&=\pi\frac{\Gamma(1)\Gamma(\alpha+\beta-2)}{\Gamma(\alpha+\beta-1)}, \Re(\alpha+\beta)>2.\nonumber
\end{align}Now, the integral over $x_3,v,z_1,z_2$ give
\begin{align}
\int_0^1x_3^{\alpha-1}(1-x_3)^{\beta-1}{\rm d}x_3&=\frac{\Gamma(\alpha)\Gamma(\beta)}{\Gamma(\alpha+\beta)},\Re(\alpha)>0,\Re(\beta)>0\nonumber\\
\int_0^1v^{\alpha-2}(1-v)^{\beta-2}{\rm d}v&=\frac{\Gamma(\alpha-1)\Gamma(\beta-1)}{\Gamma(\alpha+\beta-2)},\Re(\alpha)>1,\Re(\beta)>1\nonumber\\
\int_{z_1}\int_{z_2}b^{\alpha+\beta-3}{\rm d}z_1\wedge{\rm d}z_2&=\pi\frac{\Gamma(\alpha+\beta-2)}{\Gamma(\alpha+\beta-1)}.\nonumber
\end{align}Then the product gives
$$\pi\frac{\Gamma(\alpha)\Gamma(\alpha-1)\Gamma(\beta)\Gamma(\beta-1)}{\Gamma(\alpha+\beta)\Gamma(\alpha+\beta-1)}
=\frac{\tilde{\Gamma_2}(\alpha)\tilde{\Gamma_2}(\beta)}{\tilde{\Gamma_2}(\alpha+\beta)}=\tilde{B}_2(\alpha,\beta)
$$for $\Re(\alpha)>1,\Re(\beta)>1$, and hence the result is verified. Explicit evaluations for $p\ge 3$ will be difficult and it is simpler to use matrix methods directly starting from $p\ge 2$.

\vskip.3cm\noindent
3.\hskip.3cm{\bf Some Special Cases of Fractional Integrals in the Complex Matrix-variate Case}

\vskip.3cm Let $f(\tilde{V})=|{\rm det}(\tilde{V})|^{-\gamma}$. Consider Kober operator of the second kind when $f(\tilde{V})$ is as given above.
\begin{align}
K_{2,\tilde{U}}^{-\alpha}f&=\frac{|{\rm det}(\tilde{U})|^{\beta}}{\tilde{\Gamma_p}(\alpha)}\int_{\tilde{V}>\tilde{U}}|{\rm det}(\tilde{V})|^{-\beta-\alpha}|{\rm det}(\tilde{V}-\tilde{U})|^{\alpha-p}f(\tilde{V}){\rm d}\tilde{V}\nonumber\\
&=\frac{|{\rm det}(\tilde{U})|^{\beta}}{\tilde{\Gamma_p}(\alpha)}\int_{\tilde{V}>\tilde{U}}|{\rm det}(\tilde{V})|^{-\beta-\alpha-\gamma}|{\rm det}(\tilde{V}-\tilde{U})|^{\alpha-p}{\rm d}\tilde{V}\nonumber\\
&=\frac{|{\rm det}(\tilde{U})|^{\beta}}{\tilde{\Gamma_p}(\alpha)}\int_{\tilde{W}>O}|{\rm det}(\tilde{W})|^{\alpha-p}|{\rm det}(\tilde{W}+\tilde{U})|^{-\beta-\alpha-\gamma}{\rm d}\tilde{W},~\tilde{W}=\tilde{V}-\tilde{U}\nonumber\\
&=\frac{|{\rm det}(\tilde{U})|^{-\alpha-\gamma}}{\tilde{\Gamma_p}(\alpha)}\int_{\tilde{W}>O}|{\rm det}(\tilde{W})|^{\alpha-p}|{\rm det}(I+\tilde{U}^{-\frac{1}{2}}\tilde{W}\tilde{U}^{-\frac{1}{2}})|^{-\beta-\alpha-\gamma}{\rm d}\tilde{W}\nonumber\\
&=\frac{|{\rm det}(\tilde{U})|^{-\gamma}}{\tilde{\Gamma_p}(\alpha)}\int_{\tilde{T}>O}|{\rm det}(\tilde{T})|^{\alpha-p}|{\rm det}(I+\tilde{T})|^{-\beta-\alpha-\gamma}{\rm d}\tilde{T},~\tilde{T}=\tilde{U}^{-\frac{1}{2}}\tilde{W}\tilde{U}^{-\frac{1}{2}}\nonumber\\
&=|{\rm det}(\tilde{U})|^{-\gamma}\frac{\tilde{\Gamma_p}(\beta+\gamma)}{\tilde{\Gamma_p}(\alpha+\beta+\gamma)},~\Re(\alpha)>p-1,~\Re(\beta+\gamma)>p-1.\nonumber
\end{align}
\vskip.3cm Consider Kober operator of the first kind.
$$K_{1,\tilde{U}}^{-\alpha}f=\frac{|{\rm det}(\tilde{U})|^{-\alpha-\beta}}{\tilde{\Gamma_p}(\alpha)}\int_{\tilde{V}<\tilde{U}}|{\rm det}(\tilde{V})|^{\beta}|{\rm det}(\tilde{U}-\tilde{V})|^{\alpha-p}f(\tilde{V}){\rm d}\tilde{V}.$$

\vskip.3cm\noindent{\bf Special case 1:}\hskip.3cm $f(\tilde{V})=|{\rm det}(\tilde{V})|^{\gamma}$. Then
\begin{align}
K_{1,\tilde{U}}^{-\alpha}f&=\frac{|{\rm det}(\tilde{U})|^{-\alpha-\beta}}{\tilde{\Gamma_p}(\alpha)}\int_{\tilde{V}<\tilde{U}}|{\rm det}(\tilde{V})|^{\beta+\gamma}|{\rm det}(\tilde{U})|^{\alpha-p}|{\rm det}(I-\tilde{U}^{-\frac{1}{2}}\tilde{V}\tilde{U}^{-\frac{1}{2}})|^{\alpha-p}{\rm d}\tilde{V}\nonumber\\
&=\frac{|{\rm det}(\tilde{U})|^{\gamma}}{\tilde{\Gamma_p}(\alpha)}\int_{O<\tilde{W}<I}|{\rm det}(\tilde{W})|^{\beta+\gamma}|{\rm det}(I-\tilde{W})|^{\alpha-p}{\rm d}\tilde{W},~\tilde{W}=\tilde{U}^{-\frac{1}{2}}\tilde{V}\tilde{U}^{-\frac{1}{2}}\nonumber\\
&=|{\rm det}(\tilde{U})|^{\gamma}\frac{\tilde{\Gamma_p}(\beta+\gamma+p)}{\tilde{\Gamma_p}(\alpha+\beta+\gamma+p)}\nonumber
\end{align}for $\Re(\alpha)>p-1,\Re(\beta+\gamma)>-1$, by evaluating by using a type-1 beta integral.

\vskip.3cm\noindent{\bf Special case 2:}\hskip.3cm $f(\tilde{V})=|{\rm det}(I-\tilde{V})|^{\gamma}$.
\begin{align}
K_{1,\tilde{U}}^{-\alpha}f&=\frac{|{\rm det}(\tilde{U})|^{-\alpha-\beta}}{\tilde{\Gamma_p}(\alpha)}\int_{\tilde{V}<\tilde{U}}|{\rm det}(\tilde{V})|^{\beta}|{\rm det}(I-\tilde{V})|^{-\gamma}|{\rm det}(\tilde{U}-\tilde{V})|^{\alpha-p}{\rm d}\tilde{V}\nonumber\\
&=\frac{1}{\tilde{\Gamma_p}(\alpha)}\int_{O<\tilde{W}<I}|{\rm det}(\tilde{W})|^{\beta}|{\rm det}(I-\tilde{W})|^{\alpha-p}|{\rm det}(I-\tilde{U}^{-\frac{1}{2}}\tilde{W}\tilde{U}^{-\frac{1}{2}})|^{-\gamma}{\rm d}\tilde{W},\nonumber
\end{align}for $\tilde{W}=\tilde{U}^{-\frac{1}{2}}\tilde{V}\tilde{U}^{-\frac{1}{2}}$. Now, evaluating by using Example 6.4 of [5] we have
$$K_{1,\tilde{U}}^{-\alpha}f=\frac{\tilde{\Gamma_p}(\beta+p)}{\tilde{\Gamma_p}(\alpha+\beta+p)}{_2F_1}(\beta+p,
\gamma;\alpha+\beta+p;\tilde{U}),O<\tilde{U}<I
$$for $\Re(\alpha)>p-1,\Re(\beta)>-1$.

\vskip.3cm\noindent{\bf Special case 3:}\hskip.3cm $f(\tilde{V})=|{\rm det}(\tilde{V})|^{\gamma}|{\rm det}(I-\tilde{V})|^{-\delta}$. Then, following through the steps in special cases 1 and 2 we have
$$K_{1,\tilde{U}}^{-\alpha}f=|{\rm det}(\tilde{U})|^{\gamma}\frac{\tilde{\Gamma_p}(\beta+\gamma+p)}{\tilde{\Gamma_p}(\alpha+\beta+\gamma+p)}
{_2F_1}(\beta+\gamma+p,\delta;\alpha+\beta+\gamma+p;\tilde{U})
$$for $O<\tilde{U}<I, \Re(\alpha)>p-1,\Re(\beta+\gamma)>-1$.

\vskip.3cm\noindent{\bf Special case 4:}\hskip.3cm We can have a hypergeometric series for $f(\tilde{V})$. For the meaning of the symbol $[a]_K$, zonal polynomial $\tilde{C}_K(\tilde{V})$ and the partition $K$, see for example [5],[13]. Let
\begin{align}
f(\tilde{V})&={_rF_s}(a_1,...,a_r;b_1,...,b_s;\tilde{V})\nonumber\\
&=\sum_{k=0}^{\infty}\sum_{K}\frac{[a_1]_K...[a_r]_K}{[b_1]_K...[b_s]_K}\frac{\tilde{C}_K(\tilde{V})}{k!}\nonumber
\end{align}for  $s\ge r$ or $r=s+1$ and $\Vert\tilde{V}\Vert <1$ where $\Vert(\cdot)\Vert $ denotes a norm of $(\cdot)$. Then

\begin{align}
K_{1,\tilde{U}}^{-\alpha}f&=\sum_{k=0}^{\infty}\sum_K\frac{[a_1]_K...[a_r]_K}{[b_1]_K...[b_s]_K k!}\frac{|{\rm det}(\tilde{U})|^{-\alpha-\beta}}{\tilde{\Gamma_p}(\alpha)}\nonumber\\
&\times \int_{\tilde{V}<\tilde{U}}|{\rm det}(\tilde{V})|^{\beta}|{\rm det}(\tilde{U}-\tilde{V})|^{\alpha-p}C_K(\tilde{V}){\rm d}\tilde{V}.\nonumber
\end{align}Take out $\tilde{U}$ from $|{\rm det}(\tilde{U}-\tilde{V})|^{\alpha-p}$ and make the transformation $\tilde{W}=\tilde{U}^{-\frac{1}{2}}\tilde{V}\tilde{U}^{-\frac{1}{2}}$. Then the integral part becomes
\begin{align}
\frac{1}{\tilde{\Gamma_p}(\alpha)}&\int_{O<\tilde{W}<I}|{\rm det}(\tilde{W})|^{\beta}|{\rm det}(I-\tilde{W})|^{\alpha-p}\tilde{C}_K(\tilde{U}^{\frac{1}{2}}\tilde{W}\tilde{U}^{\frac{1}{2}})\nonumber\\
&=\frac{\tilde{\Gamma_p}(\beta+p)}{\tilde{\Gamma_p}(\alpha+\beta+p)}\frac{[\beta+p]_K}{[\alpha+\beta+p]_K}
\tilde{C}_K(\tilde{U})\nonumber
\end{align}by using (6.1.21) of [5]. Now, writing the result as a hypergeometric function we have
$$K_{1,\tilde{U}}^{-\alpha}f=\frac{\tilde{\Gamma_p}(\beta+p)}{\tilde{\Gamma_p}(\alpha+\beta+p)}
{_{r+1}F_{s+1}}(a_1,...,a_r,\beta+p;b_1,...,b_s,\alpha+\beta+p;\tilde{U})
$$for $s\ge r$ or $r=s+1$ and $\Vert\tilde{U}\Vert <1$.

\vskip.3cm\noindent{\bf Acknowledgement}
\vskip.3cm The author would like to thank the Department of Science and Technology, Government of India for the financial assistance for this work under Project No SR/S4/MS:287/05 and the Centre for Mathematical and Statistical  Sciences India for the facilities.

\vskip.3cm\noindent\begin{center}
{\bf References}
\end{center}
\vskip.2cm\noindent[1]~~R. Gorenflo, F. Mainardi, Fractional
 calculus integral and differential equations of fractional order,
 in A. Carpinteri and F. Mainardi (editors), Fractal and Fractional
 Calculus in Continuum Mechanics, Wien and New York, Springer-Verlag,
 1997, pp.223-276.

 \vskip.2cm\noindent[2]~~R. Gorenflo, F.Mainardi, Approximation of
 L\'evy-Feller diffusion by random walk, Journal for Analysis and its
 Applications, 18(2)(1999) 1-16.

\vskip.2cm\noindent[3]~~Hans J. Haubold, A.M. Mathai, R.K. Saxena,
Solutions of certain fractional kinetic equations and a fractional
diffusion equation, Journal of Mathematical Physics, 51(2010)
103506-1-103506-8.

\vskip.2cm\noindent[4]~~S. Jespersen, R. Metzler, H.C. Forgedby,
 L\'evy flights in external force fields: Langevin and fractional
 Fokker-Planck equations and their solutions, Physical Review E,
 59(3)(1999) 2736-2745.

 \vskip.2cm\noindent[5]~~A.M. Mathai, Jacobians of Matrix
 Transformations and Functions of Matrix Argument, World Scientific Publishing,
 New York, 1997.

 \vskip.2cm\noindent [6]~~A.M. Mathai, A pathway to matrix variate
 gamma and normal densities, Linear Algebra and its
 Applications, 396(2005) 317-328.

 \vskip.2cm\noindent[7]~~A.M. Mathai, Some properties of
 Mittag-Leffler functions and matrix-variate analogues: A
 statistical perspective, Fractional Calculus \& Applied Analysis,
 13(1)(2010) 113-132.

 \vskip.2cm\noindent [8]~~A.M. Mathai, Fractional integral operators in the complex matrix variate case, Linear Algebra and Its Applications, 439 (2013), 2901-2913.

 \vskip.2cm\noindent [9]~~ A.M. Mathai, Fractional integral operators involving many matrix variables, Linear Algebra and Its Applications, 446 (2014), 196-215.

 \vskip.2cm\noindent [10]~~A.M. Mathai, Hans J. Haubold, Pathway
 model, superstatistics, Tsallis statistics and a generalized
 measure of entropy, Physica A, 375(2007) 110-122.

 \vskip.2cm\noindent[11]~~A.M. Mathai, Hans J. Haubold, Special
 Functions for Applied Scientists, Springer, New York, 2008.

 \vskip.2cm\noindent[12]~~A.M.Mathai, Hans J. Haubold, Fractional
 operators in the matrix variate case, Fractional Calculus \&
 Applied Analysis, 16(2)(2013) 469-478.

 \vskip.2cm\noindent[13]~~A.M. Mathai, S.B. Provost, T. Hayakawa,
 Bilinear Forms and Zonal Polynomials, Springer, New York, [Lecture
 Notes] 1995.

 \vskip.2cm\noindent[14]~~K.S. Miller, B. Ross, An Introduction to
 the Fractional Calculus and Fractional Differential Equations,
 Wiley New York, 1993.

 \end{document}